\documentclass[12pt]{article}
\usepackage{amssymb,amsfonts,amsmath,latexsym,epsf,tikz,url}

\newtheorem{theorem}{Theorem}[section]

\newtheorem{corollary}[theorem]{Corollary}

\newcommand{\proof}{\noindent{\bf Proof.\ }}
\newcommand{\qed}{\hfill $\square$\medskip}

\begin{document}

\title{The $k$-independent graph of a graph}

\author{Davood Fatehi$^{a}$,
Saeid Alikhani  $^{a,}$\footnote{Corresponding author}, 
Abdul Jalil M. Khalaf$^b$}

\date{\today}

\maketitle

\begin{center}
$^a$Department of Mathematics, Yazd University, 89195-741, Yazd, Iran\\
{\tt  davidfatehi@yahoo.com, alikhani@yazd.ac.ir}
\medskip

$^b$Department of Mathematics,
Faculty of Computer Science and Mathematics\\ University of Kufa,
PO Box 21, Najaf, Iraq\\
{\tt abduljaleel.khalaf@uokufa.edu.iq}

\end{center}

\begin{abstract}
Let $G=(V,E)$ be a simple graph.  A set $I\subseteq V$ is an  independent set, if
no two of its members are adjacent in  $G$. The $k$-independent graph of $G$, $I_k (G)$, is defined to be the graph whose vertices correspond to the independent sets of $G$ that have cardinality at most $k$. Two vertices in $I_k(G)$ are adjacent if and only if the corresponding independent sets of $G$ differ by either adding or deleting a single vertex. In this paper, we obtain some properties of $I_k(G)$ and compute it for
some  graphs.
\end{abstract}

\noindent{\bf Keywords:}  independence number; $k$-independent graph; reconfiguration.

\medskip
\noindent{\bf AMS Subj.\ Class.:} 05C60, 05C69

\section{Introduction}
Given a simple graph $G=(V,E)$, a set $I\subseteq V$ is an independent set of $G$, if there is no edge of $G$ between any two vertices of $I$.
A maximal independent set is an independent set that is not a proper subset of any other independent set.
A maximum independent set is an independent set of greatest cardinality for  $G$. This cardinality is called independence number of $G$,
and is denoted by $\alpha (G)$.
 Reconfiguration problems have been studied often in recent years. These arise in settings
where the goal is to transform feasible solutions to a problem in a step-by-step manner, while
maintaining a feasible solution throughout.

 For the study of dominating  set reconfiguration problem: given two dominating sets $S$ and $T$ of a graph $G$, both of size at most $k$, is it possible to transform $S$ into $T$ by adding and removing vertices one-by-one, while maintaining a
dominating  set of size at most $k$ throughout?
 Regarding to this dominating set reconfiguration problem, recently the $k$-dominating graph of a graph $G$ has defined in \cite{Hass}.
The $k$-dominating graph of $G$, $D_k (G)$, is defined to be the graph whose vertices correspond to the dominating sets of $G$
that have cardinality at most $k$. Two vertices in $D_k(G)$ are adjacent if and only if the corresponding dominating
sets of $G$ differ by either adding or deleting a single vertex.  Authors in \cite{Hass}, gave conditions that
ensure $D_k(G)$ is connected. In \cite{davood} authors proved  that if $G$ is a graph without isolated vertices of order $n\ge 2$ and with $G\cong D_k(G)$, then $k=2$ and $G=K_{1,n-1}$ for some $n\ge 4$. It is also proved that for a given $r$ there exist only a finite number of $r$-regular, connected dominating graphs of connected graphs (\cite{davood}).

  One of the most well-studied problem in reconfiguration problems, is the reconfiguration of independent
sets. For a graph $G$ and integer $k$, the independent sets of size at least/exactly $k$ of $G$
form the feasible solutions. Independent sets are also called token configurations, where the independent set vertices are viewed as tokens  \cite{Bonsama}. Deciding for existence of  a reconfiguration between two
$k$-independent sets with at most $\ell$ operations is strongly NP-complete (\cite{Kami}).
 Bonamy and  Bousquet in \cite{Bonamy} have considered the  $k$-TAR reconfiguration graph, $TAR_k(G)$, as follows:

 A $k$-independent set of $G$ is a set $S\subseteq V$ with $|S|\geq k$, such that no two elements of $S$ are adjacent.
Two $k$-independent sets $I$ and $J$ are adjacent if they differ on exactly one vertex.
This model is called the Token Addition and Removal (TAR). Authors in \cite{Bonamy} provided a cubic-time algorithm to
decide whether $TAR_k(G)$ is connected when $G$ is a graph which does not contain induced paths of length $4$. Their work solves an open question
in \cite{Bonsama}. Also they  described a linear-time algorithm which decides
whether two elements of $TAR_k(G)$ are in the same connected component. As usual we denote the complete graph, path and cycle of order $n$ by $K_n$, $P_n$ and $C_n$, respectively. Also
$K_{1,n}$ is the star graph with $n+1$ vertices. 

\medskip

In the next section, we  study  the $k$-independent graph of a graph $G$. In Section 3, we study the $\alpha$-independent graph of a graph. Finally in Section 3, we exclude the empty set from the family set of independent sets of $G$, denote the new $k$-independent graph of $G$ by $I_k^*(G)$ and study its connectedness.  

\section{The $k$-independent graph of a graph}

\begin{figure}[h]
	\begin{center}
		\begin{minipage}{5cm}
			\hspace{4cm}
			\includegraphics[width=3cm,height=5cm]{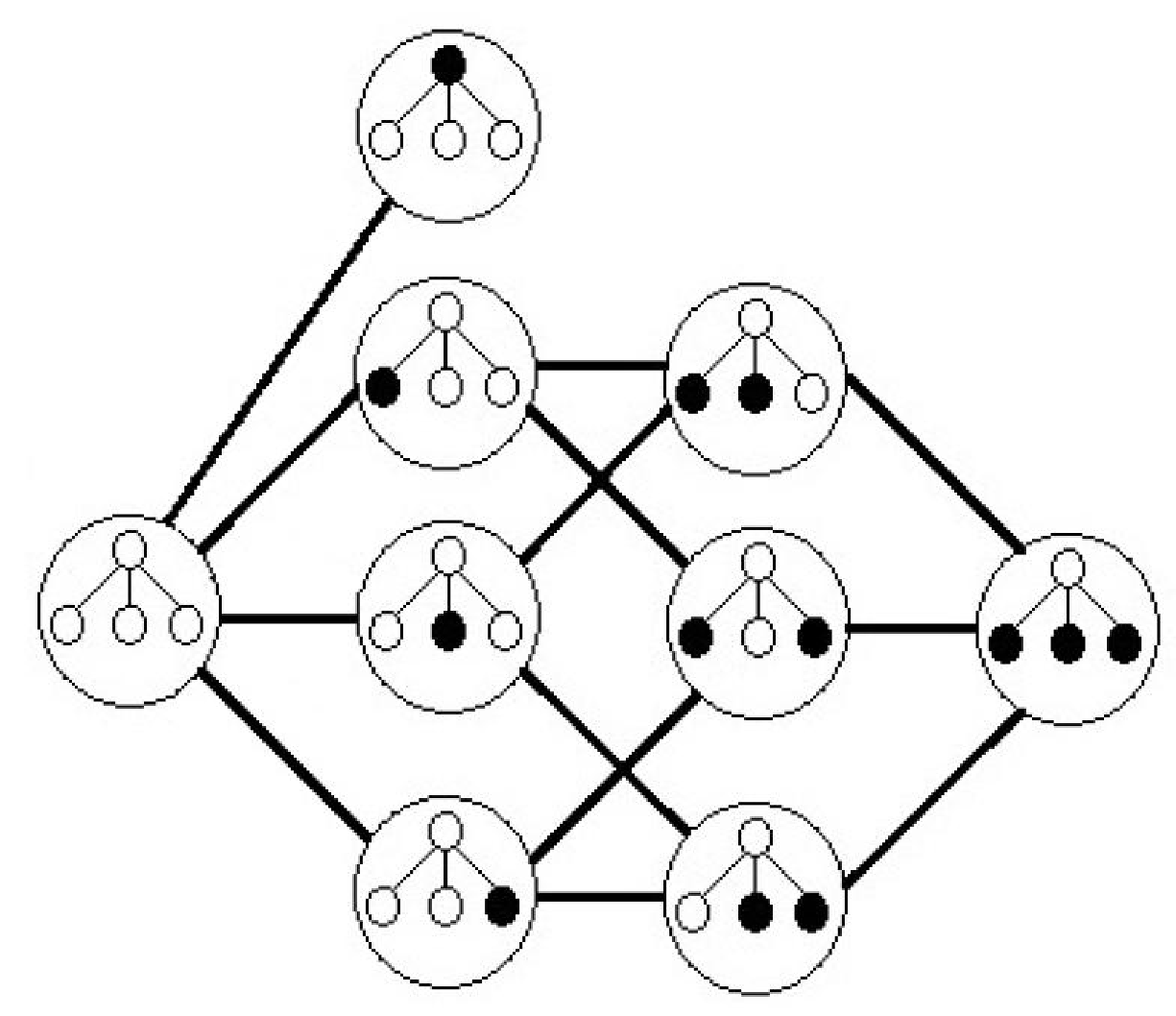}
		\end{minipage}
		\begin{minipage}{5cm}
			\includegraphics[width=3.5cm,height=5cm]{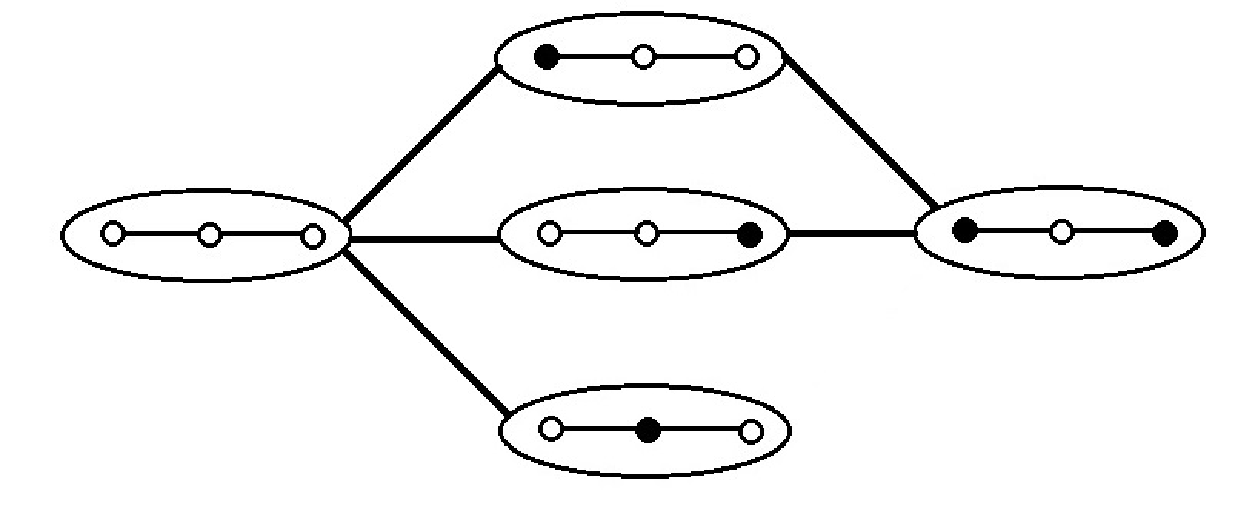}
		\end{minipage}
		\caption{ \label{figure1} Graphs $I_{3}(K_{1,3})$ and  $I_{2}(P_{3})$, respectively.}
	\end{center}\end{figure}

In this section we shall study the $k$-independent graph of a graph $G$. First let to rewrite the definition of the reconfiguration graph $TAR_k(G)$, as follows.
 For  a graph $G$ and a non-negative integer $k$, the $k$-independent graph of $G$, $I_k (G)$, is defined to be the graph whose vertices correspond to the independent  sets of $G$ that have cardinality
at most $k$. Two vertices in $I_k(G)$ are adjacent if and only if the corresponding independent
sets of $G$ differ by either adding or deleting a single vertex.
 As an example, Figure \ref{figure1} shows  $I_{3}(K_{1,3})$.

   Note that $k$-dominating and $k$-independent graph are similar to recent work in graph colouring, too. Given a graph
$H$ and a positive integer $k$, the $k$-colouring graph of $H$, denoted $G_k(H)$, has vertices
corresponding to the (proper) $k$-vertex-colourings of $H$. Two vertices in $G_k(H)$ are
adjacent if and only if the corresponding vertex colourings of $G$ differ on precisely one
vertex. Authors in \cite{4,5,6,8} studied the connectedness of $k$-colouring graphs. Also they studied their hamiltonicity. Let to introduce a notation. Let  $A$ and $B$ be independent sets of $G$ of cardinality at most $k$. We use the notation $A\leftrightarrow B$, if there is a path in $I_{k} (G)$ joining $A$ and $B$.
It is easy to see that for every $A, B \in I_{k} (G)$, 	$A \leftrightarrow B$ if and only if $B \leftrightarrow A$ and if $A \supseteq B$, then $A \leftrightarrow B$ and $B \leftrightarrow A$.

The following theorem, gives some properties of the $k$-independent graph of a graph:

\begin{theorem} \label{pro}
	\begin{enumerate}
		\item[(i)]
		If $G$ is a graph of order $n$, then $I_1(G)\cong K_{1,n}$.
		\item[(ii)]
		For every graph $G$ and every $0\leq k\leq \alpha(G)$, the independent graph $I_k(G)$ is connected and $\Delta(I_k(G))=|V(G)|$.
		\item[(iii)]
		For every graph $G$, the independent graph $I_k(G)$ is a bipartite graph.
		\item[(iv)]
		If $G\ncong \overline{K_n}$, then $I_k(G)$ is not a regular graph.
		\item[(v)]
		If $G\ncong \overline{K_n}$ then $I_k(G)$ is not a vertex-transitive graph, and so is not a Cayley graph.
	\end{enumerate}
\end{theorem}
\proof
\begin{enumerate}
	\item[(i)] It follows from the definition. 
	 \item[(ii)] It is straightforward.  
	 
	\item[(iii)]
	Let $X$ be the set of independent sets of size less than $k+1$ of $G$ with odd cardinality and $Y$ be the set of independent sets of  size
	less than $k+1$ with even cardinality. It is  clear that $X\cup Y=V(I_k(G))$ and $X\cap Y=\phi$. Suppose that $A,B\in X$,
	then $(A\backslash B)\cup (B\backslash A)$ cannot be a vertex of $I_k(G)$.  Because $|A|=|B|$ or $\big||A|-|B|\big|\geq 2$.
	So $AB$ is not an edge of $I_k(G)$ and with similar argument we have this for two vertices in $Y$.
	Therefore $I_k(G)$ is a bipartite graph with parts $X$ and $Y$.
	\item[(iv)]
	Let $G$ be a graph of order $n$. The empty set is an independent set of $G$ which has degree $n$ in $I_k(G)$.
	Let $I_1$ be an independent set of $G$ with $|I_1|=\alpha(G)$. We know that $I_1$ is adjacent to $\alpha$ independent sets.
	Since $G\ncong \overline{K_n}$, we have  $\alpha(G)\neq n$. Therefore $I_k(G)$ is not a regular graph.
	\item[(v)]
	It follows from Part (iv).\quad\qed
\end{enumerate}

\begin{theorem}  \label{tree}
	\begin{enumerate}
		\item[(i)]
		Let $G$ be a graph of order $n$. There is no integer $k$, such that $I_k(G)\cong G$.
		\item[(ii)]
		If $G\ncong K_n$, then the girth of $I_k(G)$ is $4$.
		\item[(iii)]
		\ Let $G\neq K_n$ be a graph. Then for all integers $k\geq 2$, $I_k(G)$ is not a tree.
	\end{enumerate}
\end{theorem}
\proof
\begin{enumerate}
	\item[(i)] Since for every integer number $k\geq 1$, $|V(I_k(G))|\geq n+1$, so we have the result.
	\item[(ii)] Let $v_1$ and $v_2$ be two non-adjacent vertices of graph $G$. So $\{v_1\}$ and $\{v_2\}$ are two independent sets of $G$ and
	therefore two vertices of $I_k(G)$. Now $\emptyset$, $\{v_1\}$, $\{v_1,v_2\}$, $\{v_2\}$, $\emptyset$ is a cycle in $I_k(G)$ and
	this is the shortest cycle in $I_k(G)$. Therefore the girth of $I_k(G)$ is 4.
	\item[(iii)] It follows from Part $(ii)$.\quad\qed
\end{enumerate}

\section{The $\alpha$-independent graph of some graphs}

 Let $G$ be a simple graph with independence number $\alpha$. Looks that in the among of  $k$-independent graph of $G$, the $\alpha$-independent
graph of $G$ is more important. In this section, we study the $\alpha$-independent graph of some graphs.
To study the $\alpha$-independent graph of $G$, we are interested to know the order of $I_{\alpha}(G)$.
Let $i_k$ be the number of independent  sets of cardinality $k$ in  $G$. The polynomial
$$I(G,x)=\sum_{k=0}^{\alpha(G)}i_kx^k,$$ is called the independence polynomial of $G$ (\cite{JAMC}).
Obviously $I(G,1)$ gives the number of all independent sets of a graph $G$.  In other words, $|V(I_{\alpha}(G))|=I(G,1)$.
 Since $I(K_n,x)=1+nx$, we have $I(K_n,1)=n+1$. Therefore we have the following easy result:

\begin{theorem}
	For any integer $k>1$, there is some connected graph $G$ such that $|V(I_{\alpha}(G))|=k$.
\end{theorem}
%
%
%

 The following theorem is about the $\alpha$-independent graph of stars:

\begin{theorem}
	\begin{enumerate}
		\item[(i)]
		The $n$-independent graph of $K_{1,n}$, i.e., $I_n(K_{1,n})$ is a bipartite graph with parts $X$ and $Y$, with $|X|=2^{n-1}$ and $|Y|=2^{n-1}+1$.
		\item[(ii)]
		The $n$-independent graph $I_n(K_{1,n})$ is not Hamiltonian.
	\end{enumerate}
\end{theorem}
\proof
\begin{enumerate}
	\item[(i)]
	Let $X$ be the set of independent sets of $K_{1,n}$ with even cardinality and $Y$ be the set of independent
	sets of odd cardinality. By Theorem \ref{pro}(iii), $I_n(K_{1,n})$ is a bipartite graph with parts $X$ and $Y$.
	Obviously $|X|=\sum_{k=0}^{\lfloor\frac{n}{2}\rfloor} {n\choose 2k}$ and  since the number of independent sets of $K_{1,n}$ is
	$I(K_{1,n},1)=2^n+1$, we have
	$|Y|=1+\sum_{k=1}^{\lfloor\frac{n}{2}\rfloor} {n\choose 2k-1}$. Therefore we have the result.
	\item[(ii)]
	Since a bipartite graph with different number of vertices in its parts is not a Hamiltonian graph,
	so the $n$-independent graph $I_n(K_{1,n})$ is not a Hamiltonian graph.\quad\qed
\end{enumerate}

 Here we consider the $\alpha$-independent of some another graphs. Figure \ref{figure1} shows the $I_2(P_3)$.

\begin{theorem}
	For every $n\in \mathbb{N}$, $\delta(I_{\alpha}(P_n))=\lfloor\frac{n}{2}\rfloor.$
\end{theorem}
\proof
 The minimum  degree of vertices of $I_{\lceil\frac{n}{2}\rceil}(P_n)$ is due to maximal independent sets of $P_n$
with minimum cardinality. These vertices are adjacent to  $n-\lceil\frac{n}{2}\rceil=\lfloor\frac{n}{2}\rfloor$ of independent
sets with less cardinality.\quad\qed

 Here we shall  obtain information on the Hamiltonicity of $\alpha$-independent of some specific graphs.
Using the value of the independence polynomial at $-1$, we have $I(G;-1)=i_0-i_1+i_2-\ldots+(-1)^{\alpha}i_\alpha=f_0(G)-f_1(G)$,
where $f_0(G)=i_0+i_2+i_4+\ldots$, $f_1(G)=i_1+i_3+i_5+\ldots$
are equal to the numbers of independent sets of even size and odd size of $G$, respectively.
$I(G,-1)$ is known as the alternating number of independent sets.  We need the following theorem:

\begin{theorem}{\rm\cite{Specific}} \label{value}
	For $n\geq 1$, the following hold:
	\begin{enumerate}
		\item[(i)] $I(P_{3n-2};-1)=0$ and $I(P_{3n-1};-1)=I(P_{3n};-1)=(-1)^n$;
		\item[(ii)] $I(C_{3n};-1)=2(-1)^n, I(C_{3n+1};-1)=(-1)^n$ and $I(C_{3n+2};-1)=(-1)^{n+1}$;
		\item[(iii)] $I(W_{3n+1};-1) = 2(-1)^n-1$ and $I(W_{3n};-1)=I(W_{3n+2};-1)=(-1)^n-1$.
	\end{enumerate}
\end{theorem}

\begin{corollary}
	For all positive integer $n$, the graphs $I_{\alpha}(P_{3n-1})$, $I_{\alpha}(P_{3n})$,  $I_{\alpha}(C_{n})$ and $I_{\alpha}(W_{n})$
	are not Hamiltonian.
\end{corollary}
\proof
We know that $I_{\alpha}(P_{n})$, $I_{\alpha}(C_{n})$ and $I_{\alpha}(W_{n})$ are bipartite graphs with parts containing the independent sets
of even and odd cardinality. By  Theorem \ref{value}, theses bipartite graphs have parts with different cardinality.
Therefore we have the result.\quad\qed

\section{Connectedness of $I_k^*(G)$}
     As we have seen in the Section 2, since the empty set is an independent set of any graph, so  the $k$-independent graph $I_k(G)$
is a connected graph. Let us to do not consider empty set in the study of $k$-independent graph.

\begin{figure}[h]
	\begin{center}
		\begin{minipage}{5cm}
			\includegraphics[width=3.6cm,height=5cm]{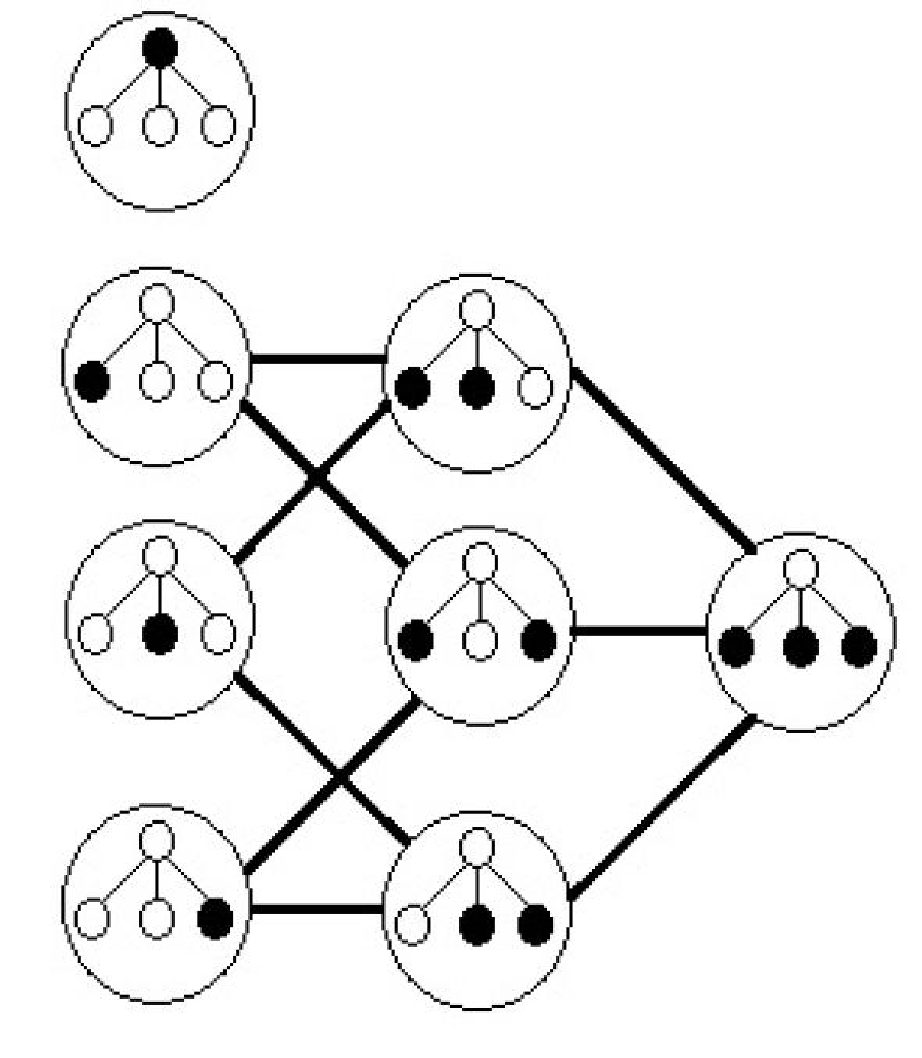}
		\end{minipage}
		\begin{minipage}{5cm}
			\includegraphics[width=3.5cm,height=5cm]{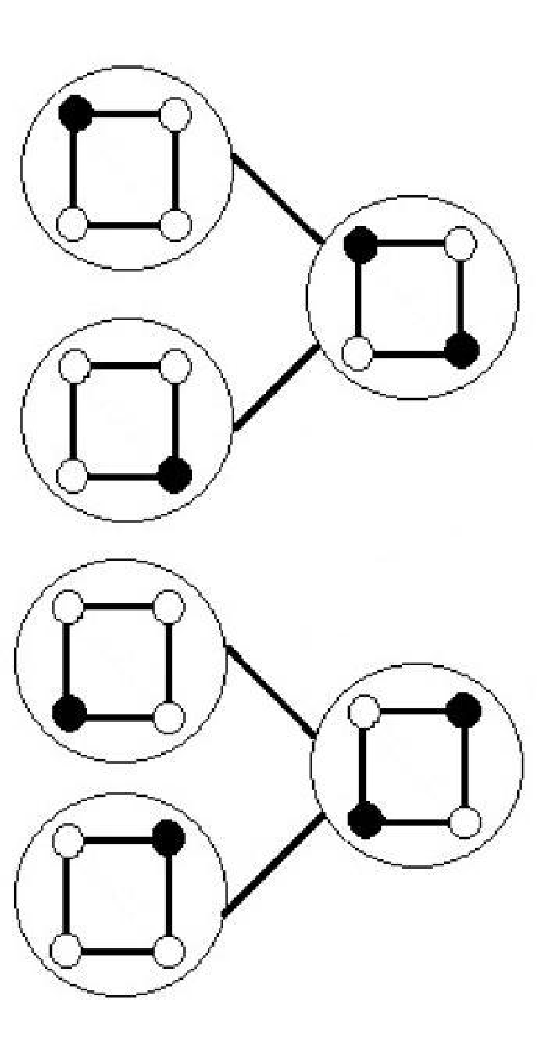}
		\end{minipage}
		\caption{\label{figure2'} Graphs  $I_{3}^*(K_{1,3})$ and $I_{2}^*(C_4)$, respectively. }
	\end{center}
\end{figure}

\begin{figure}[h]
	\begin{center}
		\begin{minipage}{5cm}
			\hspace{4cm}
			\includegraphics[width=3cm,height=5cm]{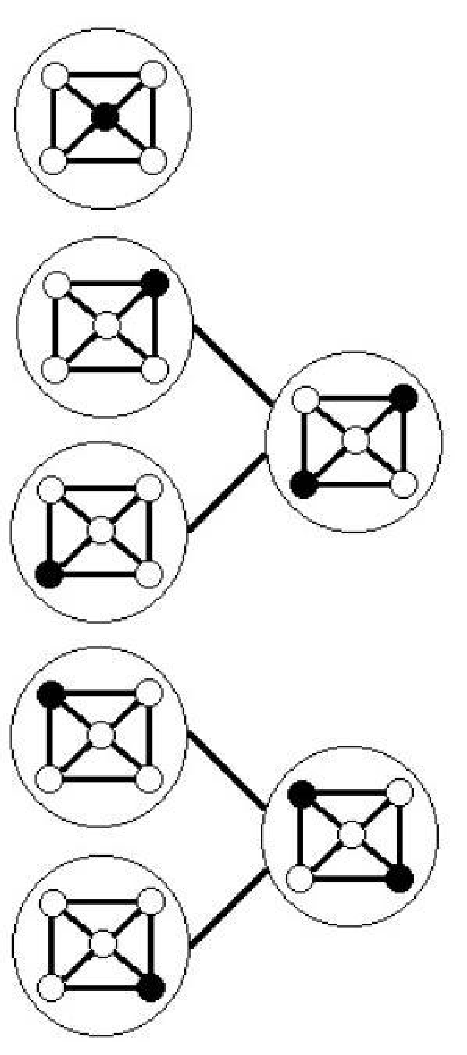}
		\end{minipage}
		\begin{minipage}{5cm}
			\includegraphics[width=3.5cm,height=5cm]{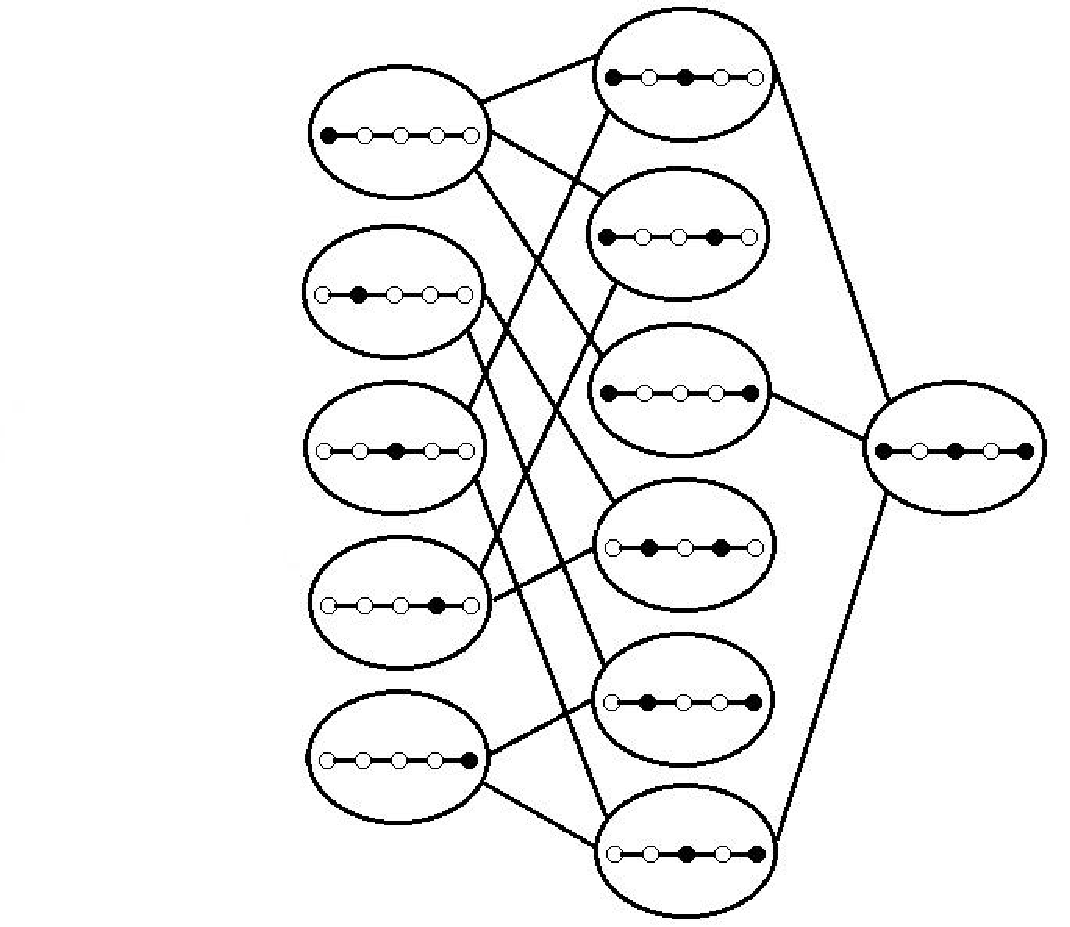}
		\end{minipage}
		\caption{\label{figure1''} Graphs $I_2^*(W_5)$ and  $I_{3}^*(P_5)$,  respectively.}
	\end{center}\end{figure}

 Suppose that $\mathcal{I}$ is a family of all independent sets of graph $G$. If we put $V(I_k(G))= \mathcal{I}\setminus \emptyset$, then
we denote the $k$-independent graph of $G$, by $I_k^*(G)$. Note that in this case,
for some $k$ and $G$,  $I_k^*(G)$  is disconnected and for some $k$ and $G$ is connected.
For example, the Figure \ref{figure2'} shows $I_3^*(K_{1,3})$ and $I_2^*(C_4)$, which are  disconnected graphs with two components.
Also Figure \ref{figure1''} shows  $I_2^*(W_5)$ and $I_{3}^*(P_5)$, respectively. Observe that $I_{3}^*(P_5)$ is connected and $I_2^*(W_5)$ is
disconnected with three components.
Theorem \ref{tree} implies that for any graph $G\neq K_n$, and  for all integers $k \geq 2$, $I_k(G)$ is not a tree, but as
we see in the Figure \ref{figure1''}, the graph $I_k^*(G)$ can be a forest. This naturally raises the question: For which graph $G$,
the component of $I_k^*(G)$ is a forest? What is the number of components?


   The following theorem is a sufficient condition for disconnectedness of $I_{\alpha}^*(G)$.
\begin{theorem}\label{dis}
	If a graph $G$ of order $n$ has a vertex of degree $n-1$, then $I_{\alpha}^*(G)$ is disconnected.
\end{theorem}
\proof
Let $v$ be a vertex of degree $n-1$. Obviously $\{v\}$ is a non-empty independent set of $G$,
and so is an isolated vertex of $I_{\alpha}^*(G)$.\quad\qed

 Note that the converse of Theorem \ref{dis} is not true. For example $I^*_2(C_4)$ has two components,
but $C_4$ is $2$-regular (Figure \ref{figure1''}).
 Now, we state the following theorem:
\begin{theorem}
	Let $K_{n_1,n_2,\ldots,n_m}$ be a complete $m$-partite  graph, then $I^*_\alpha(K_{n_1,n_2,\ldots,n_m})$ has $m$ connected components.
\end{theorem}
\proof
 Let $X_1$ and $X_2$ be two arbitrary parts of $K_{n_1,n_2,\ldots,n_m}$. Suppose that $I_1$ contains all nonempty subsets of
part $X_1$ and  $I_2$ contains all  nonempty sets of part $X_2$. Obviously, each member  of $I_1$ and each member of $I_2$ are independent sets
of $K_{n_1,n_2,\ldots,n_m}$ and so they are vertices of $I^*_\alpha(K_{n_1,n_2,\ldots,n_m})$. No member of $I_1$ is adjacent to
a member of $I_2$ in $I^*_\alpha(K_{n_1,n_2,\ldots,n_m})$. So $I^*_\alpha(K_{n_1,n_2,\ldots,n_m})$ is a disconnected graph. Since the members of
$I_1$ (and the members of $I_2$) form  a connected graph,  therefore we have $m$ components.\quad\qed

It is  obvious that, for all graph $G$ with $\alpha(G)=2$, $I^*_2(G)$ is a forest.

\begin{theorem} \label{forest}
	For all graph $G$ with $\alpha(G)>2$, the components of $I^*_k(G)$, $2\leq k\leq\alpha$, are not forest.
\end{theorem}
\proof We consider two following cases: 

 Case 1. If $k=2$. Let $\{v_1,v_2,v_3\}$ be an independent set of $G$. So $\{v_1\}$, $\{v_2\}$, $\{v_3\}$, $\{v_1,v_2\}$, $\{v_1,v_3\}$ and $\{v_2,v_3\}$ are independent sets of $G$ and vertices of $I^*_k(G)$. Therefore $\{v_1\}$, $\{v_1,v_2\}$, $\{v_2\}$, $\{v_2,v_3\}$, $\{v_3\}$, $\{v_1,v_3\}$, $\{v_1\}$ make a cycle in $I^*_k(G)$.
 
 Case 2. If $k>2$. Let $\{v_1,v_2,v_3\}$ be an independent set of $G$. So $\{v_1\}$, $\{v_1,v_2\}$ and $\{v_1,v_3\}$ are independent sets of $G$ and vertices of $I^*_k(G)$. Therefore $\{v_1\}$, $\{v_1,v_2\}$, $\{v_1,v_2,v_3\}$, $\{v_1,v_3\}$, $\{v_1\}$ make a cycle in $I^*_k(G)$ and so $I^*_k(G)$ is not a forest.\qed

Note that  if $G$ is a graph of order $n$ with  $\alpha(G)>2$, then similar to Theorem \ref{forest},  $I_k^*(G)$ cannot be a path, cycle and a chordal graph. 

\begin{theorem}
	Let $G$ be a (non complete) bipartite graph of order $n>4$. Then  $I^*_k(G)$ is connected.
\end{theorem}
\proof
Let $I_1$ and $I_2$ be two independent sets of $G$ and $|I_1|,|I_2|\leq k$, so $I_1$ and $I_2$ are two vertices of $I_k(G)$. If $I_1\cap I_2\neq \phi$ then $I_1\leftrightarrow I_1\cap I_2 \leftrightarrow I_2$.
If $I_1\cap I_2=\phi$, we consider two following cases: 

Case 1. There are $v_1\in I_1$ and $v_2\in I_2$ such that $v_1$ and $v_2$ are not adjacent then $I_1\leftrightarrow \{v_1\}\leftrightarrow \{v_1,v_2\}\leftrightarrow \{v_2\}\leftrightarrow I_2$.

Case 2. For all $v_1\in I_1$ and $v_2\in I_2$, $v_1$ is adjacent to $v_2$. So $I_1\subset A$ and $I_2\subset B$, where $A$ and $B$ are two parts of $G$.  Since $G$ is not complete bipartite graph so $I_1\neq A$ and $I_2\neq B$ and there are $v_3\in A$ and $v_4\in B$ such that $v_3\notin I_1$ and $v_3$ is not adjacent to $v_4$. We put $I_3=\left(I_1\backslash\{v_1\}\right)\cup\{v_3\}$. So $|I_3|=|I_1|$ and $I_1\leftrightarrow I_1\backslash \{v_1\}\leftrightarrow I_3$ and $I_3\leftrightarrow\{v_3\}\leftrightarrow\{v_3,v_4\}\leftrightarrow\{v_4\}\leftrightarrow I_2$. Therefore $I_1\leftrightarrow I_2$.

\end{document}